\documentclass[a4paper,reqno,12pt]{amsart}

\usepackage{amsmath, amsthm, amssymb}
\usepackage{amsfonts,mathtools}
\usepackage{bbm} 
\usepackage{epsfig,mathrsfs,graphics}
\usepackage{dsfont,graphicx,adjustbox,upgreek,wrapfig,caption,subcaption,tikz}
\usepackage{epstopdf,latexsym,colortbl,shuffle}

\usepackage[all]{xy} \entrymodifiers={!!<0pt,0.7ex>+}

\usepackage{datetime2}
\usepackage{url}
\usepackage{layouts}
\usepackage{enumitem}

\newtheorem{thm}{Theorem}

\newtheorem{lem}{Lemma}[section]

\theoremstyle{definition}
\newtheorem{rem}{Remark}[section]

\numberwithin{equation}{section}

\newcounter{parag}[subsection]

\newcounter{parage}[section]
\renewcommand{\theparage}{{\bf\thesection.\arabic{parage}}}
\newcommand{\parage}{\medskip \addtocounter{parage}{1} 
\noindent{\theparage\ } }

\newcounter{paraga}

\newcommand{\bg}{\bigskip}

\newcommand{\ga}{\gamma}
\newcommand{\Ga}{{\Gamma}}
\newcommand{\de}{\delta}

\newcommand{\De}{\Delta}
\newcommand{\eps}{{\varepsilon}}

\newcommand{\la}{\lambda}

\newcommand{\om}{\omega}

\newcommand{\Sig}{{\Sigma}}

\newcommand{\tht}{{\theta}}

\newcommand{\ph}{\varphi}

\newcommand{\ze}{{\zeta}}

\newcommand{\demi}{\frac{1}{2}}
\newcommand{\dem}{\tfrac{1}{2}}
\newcommand{\tdemi}{\frac{3}{2}}
\newcommand{\tdem}{\tfrac{3}{2}}

\newcommand{\cont}{\operatorname{cont}}

\newcommand{\I}{{\mathrm i}}
\newcommand{\dd}{{\mathrm d}}

\newcommand{\ee}{\mathrm e}

\newcommand{\ii}{^{-1}}

\newcommand{\IM}{\mathop{\Im m}\nolimits}
\newcommand{\RE}{\mathop{\Re e}\nolimits}

\newcommand{\ie}{{\emph{i.e.}}\ }

\newcommand{\eg}{{\it e.g.}\ }
\newcommand{\resp}{{resp.}\ }
\newcommand{\wrt}{{with respect to}}

\newcommand{\rhs}{{right-hand side}}

\newcommand{\dst}{\displaystyle}

\newcommand{\C}{\mathbb{C}}

\newcommand{\R}{\mathbb{R}}

\newcommand{\Z}{\mathbb{Z}}

\newcommand{\HKL}{\mathcal{H}}

\newcommand{\cL}{\mathcal{L}}

\DeclarePairedDelimiter\abs{\lvert}{\rvert}%

\newcommand{\singz}[1]{ {\operatorname{sing}}_{0}\hspace{-.07em}\big( #1 \big)}
\newcommand{\var}{\operatorname{var}}
\newcommand{\bem}{\mbox{}^\flat\hspace{-.8pt}}

\newcommand{\htb}[1]{\raisebox{-.23ex}{${\stackrel{
            \raisebox{-.18ex}{$\scriptscriptstyle\wedge$}
          }{#1}
     }$}}

\newcommand{\chb}[1]{\raisebox{-.23ex}{${\stackrel{
            \raisebox{-.23ex}{$\scriptscriptstyle\vee$}
          }{#1}
     }$}}
\newcommand{\rchb}[1]{\raisebox{-.23ex}{${\stackrel{ \hspace{-.015em}
            \adjustbox{reflect,raise=-.9ex,scale={1}{.6}}{$\upnu$}
          }{#1}
     }$}}

\newcommand{\trb}[1]{\raisebox{-.23ex}{${\stackrel{
            \raisebox{-.23ex}{$\scriptscriptstyle\triangledown$}
          }{#1}
     }$}}
\newcommand{\htn}[1]{\raisebox{.20ex}{${\stackrel{
            \raisebox{-.20ex}{$\scriptscriptstyle\wedge$}
          }{#1}
     }$}}
\newcommand{\chn}[1]{\raisebox{.20ex}{${\stackrel{
            \raisebox{-.20ex}{$\scriptscriptstyle\vee$}
          }{#1}
     }$}}
\newcommand{\trn}[1]{\raisebox{.20ex}{${\stackrel{
            \raisebox{-.20ex}{$\scriptscriptstyle\triangledown$}
          }{#1}
     }$}}
\newcommand{\rchn}[1]{\raisebox{.20ex}{${\stackrel{ \hspace{-.015em}
            \adjustbox{reflect,raise=-.9ex,scale={1}{.6}}{$\upnu$}
          }{#1}
     }$}}

\newcommand{\htg}[1]{\raisebox{-.31ex}{${\stackrel{
            \raisebox{-.23ex}{$\scriptscriptstyle\wedge$}
          }{#1}
     }$}}

\newcommand{\Dtr}[1]{{\stackrel{\raisebox{-.23ex}{$\scriptscriptstyle\triangledown$}}{#1}}}

\newcommand{\defeq}{\coloneqq} 

\newcommand{\col}{\colon\thinspace}          

\newcommand{\gB}{\mathscr B}       
\newcommand{\gD}{\mathscr D}       
\newcommand{\gL}{\mathscr L}       

\newcommand{\Lrchn}[1]{\raisebox{0ex}{${\stackrel{ \hspace{-.015em}
            \adjustbox{reflect,raise=-.9ex,scale={1}{.6}}{$\upnu$}
          }{#1}
     }$}}
\newcommand{\rchnLth}{\Lrchn{\!\cL}\hspace{.01em}^\tht}
\newcommand{\rchnLintI}{\Lrchn{\!\cL}\hspace{.01em}^\intI}

\newcommand{\eith}{\ee^{\I\tht}}

\newcommand{\eipi}{\ee^{\I\pi}}
\newcommand{\emipi}{\ee^{-\I\pi}}

\newcommand{\begla}{\begin{equation}}
\newcommand{\beglab}[1]{\begin{equation}	\label{#1}}
\newcommand{\edla}{\end{equation}}

\newcommand{\iimp}{\;\Rightarrow\;}
\newcommand{\imp}{\ens\Rightarrow\ens}
\newcommand{\Imp}{\quad\Rightarrow\quad}

\newcommand{\intI}{{(\!-\frac{\pi}{2},\frac{\pi}{2})}}
\newcommand{\tintI}{{(\!-\tfrac{\pi}{2},\tfrac{\pi}{2})}}

\newcommand{\ti}{\tilde}

\newcommand{\ens}{\enspace}

\newcommand{\ttpd}{{\tfrac{\pi}{2}}}
\newcommand{\tpd}{{\frac{\pi}{2}}}
\newcommand{\trpd}{{\frac{3\pi}{2}}}
\newcommand{\ttrpd}{{\tfrac{3\pi}{2}}}
\newcommand{\tcpd}{{\tfrac{5\pi}{2}}}

\newcommand{\Itpd}{\frac{\I\pi}{2}}

\newcommand{\Rp}{\R_{>0}}

\newcommand{\Rnp}{\R_{\le0}}

\newcommand{\datestamp}{{\small{File:\ens\hbox{\tt\jobname.tex}
\ens \DTMnow}}}

\title[Variations on the Resurgence of $\Ga$]{Variations
  on the Resurgence of \\ the Gamma Function}

\author{D.~Sauzin}
\thanks{CNRS -- 
  Observatoire de Paris, PSL Research University, 75014 Paris, and CNU
Beijing}
\thanks{{\tt David.Sauzin@obspm.fr}}
\thanks{\datestamp}

\usepackage[backend=bibtex]{biblatex}
\begin{document}

\maketitle

\begin{abstract}
We review \'Ecalle's formalism of minors, natural-majors and real-majors, and
  provide explicit formulas in the Borel plane that show the
  resurgence of the exponential of the Stirling series.
  We also discuss its Stokes phenomena in the framework of alien calculus.
\end{abstract}



\section{Introduction}


Let~$\ti\C$ denote the Riemann surface of the logarithm. The functions
%
\begla 
\la(z) \defeq \frac{\Ga(z)}{\sqrt{2\pi}\, z^{z-\demi} \ee^{-z}},
\qquad \la_c(z) \defeq z^{-c}\la(z) \quad (c\in\C)
%
%
\edla
%
%
are meromorphic in~$\ti\C$ with their principal branches holomorphic
in $\C - \Rnp$.
We have $\la(z)\xrightarrow[z\to\infty]{}1$ along~$\Rp$, with
$\log\!\big(\la(z)\big)$
asymptotic to the Stirling series 
\begla
\ti\mu(z) \defeq \sum_{n=0}^\infty \frac{B_{2n+2}}{(2n+2)(2n+1)} z^{-2n-1}.
%
%
\edla
We are interested in the asymptotic expansion of~$\la(z)$ itself, as
well as that of $\la_c(z)$ or $1/\la_c(z)$,
in the framework of \'Ecalle's Resurgence Theory.
Our study was prompted by M.~Kontsevich, who communicated us an
explicit formula for the Laplace transform of a function related
to~$\la$, which entailed resurgence~\cite{MK}.
Our purpose here is to give various formulas for
``minors'', ``natural-majors'' and ``real-majors''
of~$\la_c$, explaining along the way the meaning and interest of this
Resurgence Theory jargon, how real-majors relate to Kontsevich's
formula
and why these explicit formulas prove the resurgent character of~$\la$ or~$\la_c$
independently of that of the Stirling series.
We will also discuss the Stokes phenomenon for~$\Ga$ in the language
of \'Ecalle's alien calculus.

Resurgence Theory originated with dynamical systems but has recently
pervaded many areas of mathematics and matematical physics.

\section{Beginning of a synopsis of Resurgence Theory} \label{paragedeb}

%
\parage  \label{itemBL}
  A \emph{resurgent function} is a
  function~$\ph(z)$ that is holomorphic in an unbounded domain, in which it can be
  obtained by Borel-Laplace summation from a \emph{resurgent
    series}~$\ti\ph(z)$.

  \parage \label{itemPowers}
  A \emph{resurgent series} is a power series~$\ti\ph(z)$
  in~$z\ii$ whose Borel transform $\htb\ph(\xi)$ (the
  formal series obtained by termwise application of 
  $\gB \col z^{-c}\mapsto \xi^{c-1}/ {\Ga(c)}$) is convergent near
  $\xi=0$ and is an \emph{endlessly continuable} germ, which we call
  ``the minor''.

  \parage \label{itemEC}
  Being \emph{endlessly continuable}, roughly speaking, means that
$\htb\ph(\xi)$ analytically continues to a (possibly multivalued) holomorphic
function without any
%
natural boundary:
the analytic continuation is possible along any path except for a
discrete set of singularities along the path, that can be
circumvented, at will, to the left or to the right (with some
stipulations---see \cite{Eca85} for the most general definition, or
\cite{CNP,KS}).
%

Borel-Laplace summation
\big($\rotatebox[origin=c]{90}{\scalebox{0.85}{$\dashrightarrow$}}$\big)
\ie Laplace transform (or
some variant) composed with Borel transform, thus gives the
function~$\ph(z)$;
the formal series~$\ti\ph(z)$ (which usually is divergent
everywhere) then appears as the asymptotic expansion at infinity
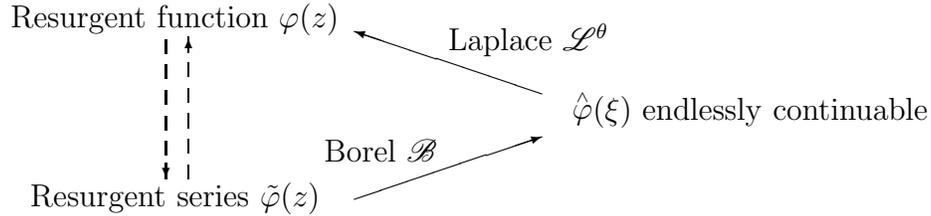
\big($\rotatebox[origin=c]{-90}{\scalebox{0.85}{$\dashrightarrow$}}$\big)
of~$\ph(z)$: 
\begin{center}
  \newlength{\DavSlength}
  \setlength{\DavSlength}{\unitlength}
  \setlength{\unitlength}{1cm}
  \begin{picture}(6,3.4)
\thinlines

\put(1.05,.4){\makebox(0,0)[br]{Resurgent series $\ti\ph(z)$}}
\put(4.4,1.8){\makebox(0,0)[l]{$\htb\ph(\xi)$ endlessly continuable}}
\put(1.3,3.2){\makebox(0,0)[tr]{Resurgent function $\ph(z)$}}

\put(1.5,.6){\vector(3,1){2.5}}
\put(4,2){\vector(-3,1){2.5}}

\put(2.6,1.1){\makebox(0,0)[br]{Borel $\gB$}}
\put(2.75,2.5){\makebox(0,0)[bl]{Laplace $\gL^\tht$}}

\multiput(-.7,.87)(0,.34){5}{\line(0,1){.185}}
\put(-.7,2.55){\vector(0,1){.2}}

\multiput(-1,1.18)(0,.34){5}{\line(0,1){.185}}
\put(-1,1.07){\vector(0,-1){.2}}
\end{picture}
  \setlength{\unitlength}{\DavSlength}
  %
  %
\end{center}
Here, we assume~$\htb\ph(\xi)$ to have analytic continuation with at
most exponential growth
 along the ray $\eith\Rp$ for some $\tht\in\R$
and 
\begin{multline}   \label{eqdefLaplth}
\ph(z) = \gL^\tht\htb\ph(z) = \int_0^{\eith\infty}
\ee^{-z\xi}\,\htb\ph(\xi)\,\dd\xi
\\[.5ex]
\text{in}\ens
\Pi^\tht_\tau \defeq \{\, z\in\ti\C \mid
-\tht-\ttpd < \arg z < -\tht + \ttpd,\ens
\RE(z\eith) > \tau
\,\}
%
%
\end{multline}
with $\tau>0$ 
large enough so that the
exponential decay of the Laplace kernel ensures the integrability of
$\ee^{-z\xi}\,\htb\ph(\xi)$ at infinity.%
\footnote{If~$\htb\ph(\xi)$ has analytic continuation with uniform
  exponential bound in a sector $\arg\xi\in I$, then
  the Laplace transforms $\gL^\tht\htb\ph(z)$, $\tht\in I$, can
  be glued together (the Cauchy theorem entail that they match) and
  give rise to a function $\gL^I\htb\ph(z)$ analytic in
  $\gD^I_\tau\defeq\bigcup\limits_{\tht\in I} \Pi^\tht_\tau$ (sectorial
  neighbourhood of~$\infty$ of opening $\abs I+\pi$).
  \label{ftnglueinggL}
  }

\parage
The most elementary examples of endlessly continuable functions are the meromorphic
functions---the analytic continuation has only one branch in this
case---and algebraic functions.
Another example is $\htn U(a,b,\xi)=\frac{\xi^{a-1}}{\Ga(a)} (1+\xi)^{b-a-1}$, whose
Laplace transform is the confluent hypergeometric function $U(a,b,z)$
%
(\cite[Ex.~6.105]{MS16}).

The Stirling series~$\ti\mu(z)$ is a resurgent series, since its Borel transform
%
\beglab{eqhatmumerom}
\htg\mu(\xi) = \sum_{n\ge0} \tfrac{B_{2n+2}}{(2n+2)!} \xi^{2n}
= \xi^{-2}\big( \tfrac{\xi}{2}\coth\tfrac{\xi}{2} - 1 \big) \in\C\{\xi\}
\edla
analytically continues as a meromorphic function on~$\C$ (with simple poles on $2\pi\I\Z^*$).
The Laplace transform $\gL^\intI\htg\mu$ is none other than $\log\!\big(\la(z)\big)$, which
is thus a resurgent function
(see \cite[pp.~244--246]{Eca85} or \cite[Th.~5.41]{MS16}).

  \parage
  We said power series in \S~\ref{paragedeb}.\ref{itemPowers}, but we may accept
non-integer powers (or even more general monomials, \eg involving powers of
$\log z$).
However, we need all monomials $z^{-c}$ of $\ti\ph(z)$ to have $\RE c>0$
for the minor $\gB\ti\ph(\xi)$
to be integrable at~$0$ and thus have a meaningful
Laplace transform.
%
%
The integrability constraint is bypassed in \'Ecalle's ``singularity
theory'' (\cite{Eca85}, \cite{Eca93}, \cite[\S3.1--3.2]{kokyu} or
\cite[Sec.~6.8]{MS16}), which not only widens the scope of Borel-Laplace
summation but also provides the appropriate framework for alien
calculus---see Section~\ref{secDirectAlien}.

A \emph{natural-major} is any endlessly continuable holomorphic
function~$\chb\ph(\xi)$ initially defined in a finite sector
$\Sig_{R,J} \defeq \{\, \xi=r\,\eith \in\ti\C \mid 0<r<R,\; \tht\in J
\,\}$, with some $R>0$ and some interval~$J$ of length $>2\pi$.
%
%
\emph{Singularities} are defined by quotienting the space of natural-majors by
the subspace $\C\{\xi\}$ of regular germs. The equivalence class
of~$\chb\ph(\xi)$ is denoted by
$\trb\ph = \singz{\chb\ph(\xi)}$.
The monodromy variation induces an operator 
\beglab{eqmonodr}
\var \col \trb\ph \mapsto \var\trb\ph=\htb\ph
\quad\text{defined by}\quad
\htb\ph(\xi) \defeq \chb\ph(\xi)-\chb\ph(\ee^{-2\pi\I}\xi).
\edla
%
%
The image~$\htb\ph$ is the \emph{minor of the singularity}~$\trb\ph$;
it is initially defined in a sector $\Sig_{R,I}$ with $I\defeq
J\cap(2\pi+J)$ but is endlessly continuable too.

A singularity~$\trb\ph$ is called an \emph{integrable singularity} if its minor~$\htb\ph$
is uniformly integrable at the origin
%
%
in every sector $\Sig_{R,I}$
and if it has a representative $\chb\ph(\xi)$
that is $o(1/\abs\xi)$ uniformly 
in every sector $\Sig_{R,I}$.
The latter condition implies that~$\trb\ph$ is uniquely determined by~$\htb\ph$;
we use the notation $\trb\ph = \bem\htb\ph$.
Conversely, any integrable
%
%
minor~$\htb\ph(\xi)$ is the minor of an integrable singularity;
\eg $\htb\ph(\xi)\in\C\{\xi\}$ \resp $\xi^{1/2}\C\{\xi\} \iimp$ a
representative of $\bem\htb\ph$ is $\chb\ph(\xi) \defeq
\htb\ph(\xi) \tfrac{\log\xi}{2\pi\I}$ \resp $\demi\htb\ph(\xi)$.

\parage
Suppose~$\trb\ph$ is an arbitrary singularity whose minor extends
analytically to a sector
$\Sig_{\infty,I}$, 
with $\abs{\htb\ph(\xi)} \le C\, \ee^{\tau\abs{\xi}}$ on
$\Sig_{\infty,I}\cap\{\abs\xi\ge1\}$.
Then, given $\tht\in I$, for any natural-major~$\chb\ph(\xi)$ and any $\de>0$ small
enough, we can consider the contour
$C_{\de,\tht} \defeq \{ \tht' \in
[\tht-2\pi,\tht] \mapsto \de\,\ee^{\I\tht'} \}$ and the function
\beglab{eqdefLaplMaj}
\cL^\tht \trb\ph(z) \defeq \int_{C_{\de,\tht}} \ee^{-z\xi}\,\chb\ph(\xi)\,\dd\xi
+ \int_{\de\,\eith}^{\eith\infty} \ee^{-z\xi}\,\htb\ph(\xi)\,\dd\xi
\quad\text{holomorphic for}\ens  z\in \Pi^\tht_\tau,
\edla
which does not depend on~$\de$ nor on the
chosen natural-major $\chb\ph(\xi)$
and can be rewritten in terms of the sole
natural-major as an integral
$\cL^\tht\chb\ph(z)\defeq\int_{\HKL_\tht} \ee^{-z\xi}\,\chb\ph(\xi)\,\dd\xi$
over a $\tht$-rotated Hankel contour $\HKL_\tht$ running from
$\ee^{\I(\tht-2\pi)}\infty$ to $\eith\infty$ and circling anticlockwise
around~$0$ at distance~$\de$
(at least if~$\chb\ph(\xi)$ itself has at most exponential growth
at~$\infty$ in directions~$\tht$ and~$\tht-2\pi$).
Moreover, as in Footnote~\ref{ftnglueinggL}, the various
$\cL^\tht \trb\ph$, $\tht\in I$, can be glued together, giving rise to
$\cL^I\trb\ph$ holomorphic in $\gD^I_\tau$.
This is an extension of the usual Laplace transform in the sense that
$\cL^I\,\bem\htb\ph = \gL^I\htb\ph$
for any integrable minor~$\htb\ph$.
%

For example, the Laplace transform of the singularity defined by
\beglab{eqdefchnIcIn}
\chn I_c(\xi) \defeq \tfrac{\xi^{c-1}}{(1-\ee^{-2\pi\I c})\Ga(c)}
\quad\text{for $c\in\C-\Z_{>0}$},
\qquad
\chn I_c(\xi) \defeq \tfrac{\xi^{c-1}}{(c-1)!} \frac{\log\xi}{2\pi\I}
\quad\text{for $c\in\Z_{>0}$}
\edla
is $\cL^I\trn I_c = z^{-c}$ for any~$I$.
In general, the asymptotic behaviour at~$\infty$ of $\cL^I\trb\ph(z)$
only depends on the asymptotic behaviour at~$0$
of a natural-major: 
if $\chb\ph(\xi)$ can be represented as a convergent series of monomials
proportional to $\chn I_{c_n}(\xi)$, then $\cL^I\trb\ph(z)$ has an
asymptotic expansion at~$\infty$ given by the corresponding series of
monomials proportional to~$z^{-c_n}$.

%
%
%
%
%

\parage
\label{subsecconvolsing}
  Resurgence and Borel summation are compatible with multiplication,
  the corresponding operation on singularities being the ``convolution
  of singularities''~$\,\Dtr *$.


\section{Main results}

We now state our main results as two theorems, which
will be proved (in Sections~\ref{secpfthmun}--\ref{secpfthmdeux}) independently one of the other and independently of the
Stirling series.

\parage
The Lambert function $x\mapsto W(x)$ is implicitly defined
by $W(x)\,\ee^{W(x)} = x$; the canonical reference for the description
of its analytic continuation is \cite{CGH96}, the only singularities
are located at $x=-\ee^{-1}$ and~$0$.
One denotes by~$W_{0}$ and~$W_{-1}$ its real branches so that
$-\ee^{-1}<x<0 \iimp W_{0}(x) > W_{-1}(x)$,
both have square root branch points at~$-\ee^{-1}$.


\begin{thm} \label{thmMinNat}
%
    %
 \textbf{\textup{(i)}}\; 
 $\la_{3/2}(z)$ is a resurgent function 
 obtained by~$\cL^\intI$ from
 the integrable singularity~$\trn\la_{3/2}$ that has natural-major and minor explicitly given,
 for $\arg\xi=0$, by 
  \beglab{eqfirstdefhtnlatdem}
  \chn\la_{3/2}(\xi) = \tfrac{1}{\sqrt{2\pi}} W_0(-\ee^{-1-\xi}),
\quad
\htn\la_{3/2}(\xi) = \tfrac{1}{\sqrt{2\pi}} \big(
W_0(-\ee^{-1-\xi})-W_{-1}(-\ee^{-1-\xi}) \big).
  \edla
 \textbf{\textup{(ii)}}\; 
  $\chi(z)\defeq\frac{1}{\la_{-3/2}(z)}$ is a resurgent
  function 
  obtained by~$\cL^\intI$ from
 the integrable singularity~$\trb\chi$ that has natural-major and minor explicitly given,
 for $\arg\xi=-\pi$, by 
  \begla
  \chb\chi(\xi) = \tfrac{\I}{\sqrt{2\pi}} W_{-1}(-\ee^{-1+\xi}),
\quad
    \htb\chi(\xi)  
%
    = \tfrac{\I}{\sqrt{2\pi}} \big( W_{-1}(-\ee^{-1+\xi}) -
    W_{0}(-\ee^{-1+\xi})\big).
    \edla
 \textbf{\textup{(iii)}}\; 
    Define inductively 
    $a_1= 1$ and $a_k= \tfrac{1}{k+1}\big(
a_{k-1} - \sum_{\ell=2}^{k-1} \ell a_\ell a_{k+1-\ell}
\big)$ for $k\ge2$, so that
\begin{gather}
\label{eqdefseqan}
  a_2 = \tfrac{1}{3},
\quad a_3 = \tfrac{1}{36},
\quad a_4 =-\tfrac{1}{270}, 
\quad a_5 =\tfrac{1}{4320}, 
\quad a_6 =\tfrac{1}{17010}, 
\quad a_7 =-\tfrac{139}{5443200}, 
\quad \ldots
\shortintertext{Then}
\label{eqexplatdem}
\;\qquad \chn\la_{3/2}(\xi) =
-\tfrac{1}{\sqrt{2\pi}}\big(1+ \sum\limits_{k\ge1} (-1)^{k} a_k(2\xi)^{k/2}\big),
%
%
\quad
\htn\la_{3/2}(\xi) =
\tfrac{1}{\sqrt{\pi}} \sum\limits_{n\ge 0} a_{2n+1} 2^{n+1} \xi^{n+\demi}
%
%
%
\shortintertext{and}
\label{eqexpchi}
\ens\qquad \chb\chi(\xi) =
-\tfrac{\I}{\sqrt{2\pi}}\big( 1 + \sum_{k\ge1} \I^k a_k (2\xi)^{k/2} \big),
\quad\ens
\htb\chi(\xi) =
\tfrac{1}{\sqrt{\pi}} \sum\limits_{n\ge 0} (-1)^n a_{2n+1} 2^{n+1} \xi^{n+\demi}.
%
%
\end{gather}
    %
%
  \end{thm}


  Notice that, as a consequence of the above statements,
  \beglab{eqsymmchilatdem}
  \chb\chi(\xi) = \I \chn\la_{3/2}(\ee^{-\I\pi} \xi),
  \qquad
  \htb\chi(\xi) = \I \htn\la_{3/2}(\ee^{-\I\pi} \xi)
  = -\I \htn\la_{3/2}(\ee^{\I\pi} \xi).
  \edla
  Since 
  the analytic continuation of the Lambert~$W$ function has no
  other singularities than those located at~$-\ee^{-1}$ and~$0$, it
  follows 
  that the minors
  $\htn\la_{3/2}(\xi), \htb\chi(\xi) \in \xi^{\demi}\C\{\xi\}$ as well
  as the natural-majors
  $\chn\la_{3/2}(\xi), \chb\chi(\xi) \in \C\{\xi^{\demi}\}$ are endlessly
  continuable; all their singularities are located above $2\pi\I\Z$
  and are square root branch points.
  On the Riemann surface of~$\xi^\demi$, we thus have four
  singularity-free sectors, corresponding to arguments in
  \[
    J_1 \defeq \tintI, \quad
    J_2 \defeq (\ttpd,\ttrpd), \quad
    J_3 \defeq (\ttrpd,\tcpd), \quad
    J_4 \defeq (\!-\ttrpd,-\ttpd).
  \]
  The functions $\la_{3/2}=z^{-\tdemi}\la$ and $\chi=\frac{z^{-\tdemi}}{\la}$
  %
  are the Laplace transforms of their minors as in
  Formula~\eqref{eqdefLaplth} with any $\tht\in J_1$ and $\tau>0$,
  hence, in the notation of Footnote~\ref{ftnglueinggL},
  \beglab{eqlatdchiasgLmin}
  z^{-\tdemi}\la(z) = \gL^{J_1}\htn\la_{3/2}(z), \quad
  \frac{z^{-\tdemi}}{\la(z)} = \gL^{J_1}\htb\chi(z)
  \quad\ens \text{in} \ens
  \bigcup_{\tht\in J_1,\,\tau>0} \Pi^\tht_\tau = \{-\pi < \arg z<\pi\}
  \edla
  (using~$J_3$ instead of~$J_1$ amounts to a simple change of branch of
  the square roots).
  The functions $\gL^{J_2}\htn\la_{3/2}$ and $\gL^{J_2}\htb\chi$ are
  defined in $\{-2\pi < \arg z<0\}$ and differ from the previous ones
  by the \emph{Stokes phenomenon} through $\tht=\tpd$, which in this
  case is equivalent to Euler's reflection formula---see Section~\ref{secDirectAlien}.

  It follows from \S~2.7 that, for every $c\in\C$,
  $\la_c = z^{-c+3/2}\la_{3/2}$ is resurgent
  with $\trn\la_c = \trn I_{c-3/2} \,\Dtr *\, \trn\la_{3/2}$. 
Note that $\trn\la_c$ is an integrable singularity if and only if $\RE c>0$.
  
  The sequence $(a_n)$ was defined in~\cite{MM90} in relation with the
  asymptotics of~$\Ga$.
  Since $\Ga(n+\frac{3}{2}) = 2^{-n-1} (2n+1)!!\, \sqrt\pi$, the
  expansions \eqref{eqexplatdem}--\eqref{eqexpchi}
   imply asymptotic expansions in the same domain as in~\eqref{eqlatdchiasgLmin}:
  \beglab{eqasymplaunsla}
  \qquad\qquad
  \la(z)\sim
\ti\la(z) \defeq \sum_{n\ge0} (2n+1)!!\, a_{2n+1} z^{-n},  
\quad
  \frac{1}{\la(z)}\sim
  \ti\la(-z) 
\quad\text{in}\ens \C-\Rnp.
  \edla
  
\parage
\emph{Real-majors}\footnote{They owe their name to their usefulness in
  the summation of divergent series with real coefficients.}
constitute a variant of the notion of natural-majors
(see \cite[pp.~81--82]{Eca93}):
a real-major of a singularity~$\trb\ph$ is a function
$\rchb\ph(\xi) = -2\pi\I\, \chb\ph(\ee^{-\pi\I}\xi)$
where $\chb\ph(\xi)$ is any natural-major. The minor and
the Laplace transform of~$\trb\ph$ can be retrieved as
\beglab{eqrchbgL}
\htb\ph(\xi) = - \tfrac{1}{2\pi\I}\big( \rchb\ph(\ee^{\pi\I}\xi) -
\rchb\ph(\ee^{-\pi\I}\xi) \big),
\quad
\cL^\tht \trb\ph(z) = \rchnLth\rchb\ph(z) \defeq
\tfrac{1}{2\pi\I} \int_{\HKL_{\tht+\pi}} \ee^{-z\xi}\,\rchb\ph(\xi)\,\dd\xi.
\edla
By Fourier-Laplace inversion, there are also integral formulas for a natural-major~$\chb\ph(\xi)$
or a real-major~$\rchb\ph(\xi)$ in terms of $\ph=\cL^\tht \trb\ph$, 
as well as for the minor~$\htb\ph(\xi)$ when~$\trb\ph$ is an integrable singularity.
For instance with $\tht=0$, in which case~$\ph$ is defined at least
in the half-plane~$\Pi^0_\tau$,
a real-major is given by
\beglab{eqrhcbphu}
\rchb\ph_{[u]}(\xi) = \int_u^{+\infty} \ee^{-z\xi}\,\ph(z)\,\dd z
\quad \text{for}\ens \arg\xi=0
\edla
for any $u\in\Pi^0_\tau$.
If $\ph(z)$ happens to be analytic along~$\Rp$ and integrable at
$z=0$, then one can take $u=0$, \ie
$\rchb\ph_{[0]}=\gL^0\ph$.


\begin{thm} \label{thmRealMaj}
  If $\RE c < \demi$, then $\la_c = \rchnLintI\rchn\la_c$
in $\C-\R_{\le0}$, with a real-major $\rchn\la_c\defeq
\gL^0\la_c$ explicitly given by 
\beglab{eqrchnlac}
\qquad\qquad\qquad\qquad
\rchn\la_c(\xi) = \tfrac{\Ga(\tdemi-c)}{\sqrt{2\pi}} 
\int_\R \, (\xi+\ee^Q-Q-1)^{c-\tdemi} \, \dd Q
\quad \text{for} \ens \xi \in \C-\R_{\le0}.
\edla
The function~$\rchn\la_c(\xi)$ extends analytically to the Riemann surface of the
logarithm except at the points that project onto $2\pi\I\Z$.
Is principal branch extends without any singularity to the sector
$\{-\trpd < \arg\xi < \trpd\}$.
%
%
\end{thm}

  
Note that for $\xi \in \C-\R_{\le0}$ the integrand in~\eqref{eqrchnlac} is well-defined, because
\beglab{ineqexp} 
\ee^Q-Q-1\ge0 \quad\text{for all}\ens Q\in\R,
\edla
and has asymptotic equivalent
$\ee^{(c-\tdemi)Q}$ as $Q\to+\infty$ and $\abs{Q}^{c-\tdemi}$ as
$Q\to-\infty$,
whence integrability follows under our condition $\RE c < \demi$.

The analytic continuation of~$\rchn\la_c(\xi)$ along any path
of~$\ti\C$ which avoids $2\pi\I\Z$ can be explicitly obtained by
continuously deforming the integration path in the \rhs\
of~\eqref{eqrchnlac}. Theorem~\ref{thmRealMaj} thus entails a
constructive proof of the resurgent character of~$\la_c(z)$,
independently of Theorem~\ref{thmMinNat}.


\begin{rem}
$\xi\notin 2\pi\I\Z$ means that
%
  %
$  \Sig_\xi \defeq \big\{ (Q,P) \in \C^2 \mid P^2 = \xi + \ee^Q - Q-
  1\big\}
$ 
  is a regular complex curve (of infinite genus).
  When 
  $c\in\Z_{\le0}$, the consideration
of~$\rchn\la_c(\xi)$ or of its monodromy around a point of $2\pi\I\Z$
leads to the periods
$\int_\ga {P^{-3+2c}}\,{\dd Q}$, $\ga \subset \Sig_\xi$.
\end{rem}


\begin{rem}
One can do everything with 
$\nu(z) \defeq \frac{\Ga(z+\demi)}{\sqrt{2\pi}\, z^{z} \, \ee^{-z}}$
and $\nu_c(z) \defeq z^{-c}\nu(z)$
instead of~$\la$ and~$\la_c$.
One obtains slightly different formulas, \eg
\begla
\RE c < 1 \imp
\gL^0\nu_c(\xi) = 2^{-\tdemi}\!\!
\int_\R \, (\xi+\ee^Q-Q-1)^{c-\tdemi} \, \ee^{Q/2} \,\dd Q
\quad \text{for} \ens \xi \in \C-\R_{\le0}.
\edla
This formula with $c=0$ is the one obtained by M.~Kontsevich~\cite{MK}.
\end{rem}


\begin{rem}
One can compute the minor
  $\htn\la_c(\xi)$ from~(\ref{eqrchbgL}\hspace{.075em}a) and~\eqref{eqrchnlac}
  using analytic continuation in~$c$. One gets nice integral formulas for~$c$
  integer,
  for instance
  \begla
  \htn\la_1(\xi) = 
  \tfrac{1}{2\pi\I\sqrt2}\oint_\ga 
{(-\xi+\ee^Q - Q-1)^{-1/2}} {\ee^Q}\, \dd Q, \quad 
\text{$\ga$ encircling~$0$.}
\edla
When~$c$ is half-integer, the link with Theorem~\ref{thmMinNat} and
the Lambert~$W$ function is that,
  for $0<\abs\xi<\!\!<1$, the two local solutions to $-\xi + \ee^Q - Q-
  1=0$ are given by 
\begla
Q_\pm(\xi) = -1-\xi-W_{0,-1}(-\ee^{-1-\xi}) = \pm(2\xi)^{1/2}
  + O(\xi).
  \edla
\end{rem}


\parage
The space of resurgent series is stable under multiplication
(\S~2.7) and nonlinear operations like
composition or substitution into a convergent series \cite{Eca85}, \cite{Sau15}, \cite{KS}.
%
%
On the other hand, nonlinear analysis is also compatible with Borel-Laplace summation, but
that is much more elementary \cite[\S~3.5]{Sau15}.
For instance, since $\log\la = \gL^{J_1} (\gB \ti\mu)$ as mentioned
earlier, it follows that
$\la = 1 + \gL^{J_1} \gB(-1+\exp{\ti\mu}) \sim \exp{\ti\mu}$.
Comparing with~\eqref{eqasymplaunsla}, by uniqueness of the asymptotic
expansion, we get $\ti\la = \exp{\ti\mu}$.

The resurgent character of~$\ti\la$ could thus be deduced from that
of the Stirling series~$\ti\mu$:
it is always the case that the exponential of a series whose Borel
transform is meromorphic, thus resurgent, has an endlessly continuable
Borel transform.
However, the resulting germ usually admits no closed form
formula, even if one started with an explicit meromorphic function:
one obtains a function that is guaranteed to be endlessly
continuable, the singularities of which can be located and analyzed by
means of \'Ecalle's alien calculus, but there is no reason why this
function should be susceptible of an expression in terms of elementary
ones.
As testified by Theorems~\ref{thmMinNat} and~\ref{thmRealMaj},
the Borel transform of the exponential of the
Stirling series is a noticeable exception.


\section{Proof of Theorem~\ref{thmMinNat}. Alien Calculus for~$\ti\la$}
\label{secMinor}
\label{secpfthmun}




\begin{lem}
For any $z\in\C$ with $\RE z>0$, 
\beglab{eqminGa}
\sqrt{2\pi}\, z^{-\frac{3}{2}} \la(z) = \int_0^{+\infty} \ee^{-z\xi}
\big(q_+(\xi)-q_-(\xi)\big)\,\dd\xi,
\edla
where $0<q_-(\xi)<q_+(\xi)$ are the roots of the equation
$q-\ln q-1=\xi$ for $\xi>0$.
\end{lem}

\begin{rem}
We noticed later that this result and its proof are
essentially in~\cite{MM90}. That paper goes on making use of the expansion
of~$q_\pm$ at $\xi=0$ but does not refer to the Lambert~$W$ function which,
however, allows one to write
\beglab{eqidentifqpmW}
q_+(\xi) = -W_{-1}(-\ee^{-1-\xi}), \quad
q_-(\xi) = -W_{0}(-\ee^{-1-\xi}).
\edla
The connection between the expansion of~$q_\pm$ at~$0$ and the Lambert
function was noticed in~\cite{BC99} and used as an illustration of ``experimental
mathematics'', as we later discovered.
\end{rem}

\begin{proof}
We start with the formula
\beglab{eqLaplMon}
z^{-c-1} = \int_0^{+\infty} \ee^{-z \ze} \frac{\ze^c}{\Ga(c+1)} \,\dd \ze
\quad \text{for $\RE c>-1$ and $\RE z>0$}
\edla
(the very basis of Borel-Laplace summation method!).
Choosing $c=z$, we get
\[
%
%
z^{-z} \Ga(z) = z^{-z-1} \Ga(z+1) = \int_0^{+\infty} \ee^{-z q} {q^z} \,\dd q
= \int_0^{+\infty} \ee^{-z (q-\ln q)} \,\dd q.
\]
The map $q\mapsto \ze = q-\ln q$ induces 
 a decreasing diffeomorphism $(0,1) \to (1,+\infty)$
and an increasing diffeomorphism $(1,+\infty) \to (1,+\infty)$.
The inverse of the decreasing diffeomorphism is
$\ze \mapsto q_-(\ze-1)$
and the inverse of the increasing one is $\ze \mapsto q_+(\ze-1)$,
both have a square root singularity at $\ze=1$.
Thus, splitting the integral in two parts and changing variable,
we get
\begin{multline*}
z^{-z}\Ga(z) = -\int_1^{+\infty} \ee^{-z\ze} q_-'(\ze-1)\,\dd\ze
\, + \int_1^{+\infty} \ee^{-z\ze} q_+'(\ze-1)\,\dd\ze \\[1ex]
= \ee^{-z}\, \int_0^{+\infty} \ee^{-z\xi} \big( q_+'(\xi)-q_-'(\xi) \big)\,\dd\xi.
\end{multline*}
We obtain~\eqref{eqminGa} by multiplying by $z\,\ee^z$ and integrating by parts.
One gets~\eqref{eqidentifqpmW} from
$q-\ln q=\ze \,\Leftrightarrow\, q\,\ee^{-q}=\ee^{-\ze}
  \,\Leftrightarrow\, -q = W_{-1}(-\ee^{-\ze})$ or
  $W_{0}(-\ee^{-\ze})$.
\end{proof}

We thus have obtained
$\la_{3/2} = \gL^0 \htn\la_{3/2}$ with~$\htn\la_{3/2}$ as in~\eqref{eqfirstdefhtnlatdem}.
As already mentioned, the endless continuability and the location of
the singularities of~$\htn\la_{3/2}$ are simple consequences of the
complex analytic structure of~$W$.
The existence of constants $A,B>0$ such that
$\abs{\htn\la_{3/2}(\xi)} \le A\abs{\xi}+B$ for $\arg\xi\in J_1\cup
J_2 \cup J_3 \cup J_4$ allows us to vary the integration direction
in the Laplace representation
of~$\la_{3/2}(z)$ and get
$\la_{3/2} = \gL^\intI \htn\la_{3/2}$.




Since $q=1$ is a critical point of multiplicity~$1$ for the
holomorphic function $q\mapsto q-\log q-1=\xi$, the critical value
$\xi=0$ is a square root branch point 
for the two local inverse branches,
%
$q_-(\xi)$ and $q_+(\xi)=q_-(\ee^{-2\pi\I}\xi)$;
these are two holomorphic functions of~$\xi^{1/2}$, one for each
choice of the square root, which admit convergent Puiseux expansions
in $\xi^{1/2}\C\{\xi^{1/2}\}$
(equivalently: $W_0$ and~$W_{-1}$ have a square root branch point at~$-1/\ee$).
Following \cite{MM90} and \cite {BC99}, we write the Puiseux expansion
in terms of $2\xi$ so as to get rational coefficients:
\begin{gather}
\label{eqPuiseuxW}
q_-(\xi) = -W_0(-\exp({-1-\xi})) = 1+ \sum_{k\ge1} (-1)^{k} a_k
(2\xi)^{k/2}, \\[1ex]
\label{eqPuiseuxWm}
q_+(\xi) = - W_{-1}(-\exp({-1-\xi})) =
-W_0(-\exp({-1-\ee^{-2\pi\I}\xi})) =
 1 + \sum_{k\ge1} a_k (2\xi)^{k/2}
\end{gather}
with the sequence $(a_k)_{k\ge1}$ of Part~(iii) of Theorem~\ref{thmMinNat}
(the induction formula is obtained in \cite{MM90} from the
differential equation satisfied by $q_\pm-1$ \wrt\ $(2\xi)^{1/2}$).
The Puiseux expansion~(\ref{eqexplatdem}\,b) for~$\htn\la_{3/2}(\xi)$
follows.




The relation $q_+(\xi)=q_-(\ee^{-2\pi\I}\xi)$ shows
that~$\htn\la_{3/2}$ is the monodromy variation of
\begla
\chn\la_{3/2}(\xi) \defeq -\frac{1}{\sqrt{2\pi}} q_-(\xi) = 
\frac{1}{\sqrt{2\pi}} W_0(-\exp({-1-\xi}))
\edla
which certainly is $o(1/\abs\xi)$ in view of~\eqref{eqPuiseuxW}
and is thus a natural-major
for~$\bem\htn\la_{3/2}$.
This concludes the proof of Part~(i) of Theorem~\ref{thmMinNat} as well as~\eqref{eqexplatdem}.

%




We will find a natural-major for $\frac{1}{\la_{-3/2}(z)}$
starting with the formula
\beglab{eqinvLaplMon}
\text{$\RE c>0$ and $\ze>0$} \Imp
\frac{\ze^c}{\Ga(c+1)} = \frac{1}{2\pi\I} \int_{u-\I\infty}^{u+\I\infty} \ee^{q \ze}\, q^{-c-1} \,\dd q
\quad \text{for any\; $u>0$}
\edla
(which can be viewed as a particular case of integral Borel transform, obtained by
Fourier-Laplace inversion from~\eqref{eqLaplMon}).
Fix $z>1$.
Choosing $\ze = z$ and $c=z-1$:
\beglab{eqintermed}
\frac{z^{z-1}}{\Ga(z)} = \frac{1}{2\pi\I}
\int_{u-\I\infty}^{u+\I\infty} \ee^{q z}\, q^{-z} \,\dd q
= \frac{1}{2\pi\I} \int_{u-\I\infty}^{u+\I\infty} \ee^{- z (-q+\log q)}\, \dd q
\quad \text{for any\; $u>0$.}
\edla
Let us fix $u>1$ and perform the change of
variable
%
%
$\ze = -q+\log q$
%
%
along the integration contour $\{ q = u + \I t \mid t \in \R \}$.
The new contour is
%
%
$\ga_u \defeq \{ \ze(t)=\ze_1(t) + \I \ze_2(t) \mid t \in \R \}$
with
\begla
\ze_1(t) \defeq -u+\dem\ln(u^2+ t^2),
\ens
\ze_2(t) \defeq -t + \arctan\tfrac{t}{u}.
\edla
We have $-q\,\ee^{-q} = - \ee^\ze$, hence
$- q = W(- \ee^\ze)$ is a branch of the Lambert $W$ function.

The even function $t\mapsto\ze_1(t)$ is decreasing on $(-\infty,0]$,
while the odd function $t\mapsto\ze_2(t)$ is decreasing;
the contour~$\ga_u$ comes from infinity in the first quadrant,
reaches the real point $\ze(0)=-u+\ln u < - 1$,
and goes to infinity in the fourth quadrant.
Since $-\ee^{\ze(0)}>-\ee^{-1}$ and $-u<-1$, we see that $-u =
W_{-1}(-\ee^{\ze(0)})$, hence
the inverse change of variable is given by
$-q = - (u+\I t) = W_{-1}(- \ee^{\ze(t)})$, at least for~$t$ close to~$0$.

One sees that, by taking $u$ close enough to~$1$, we can ensure that the contour~$\ga_u$
crosses the line $\{\IM\ze=-1\}$ the first time between $-1$ and
$-1+2\pi\I$ and the second time between $-1-2\pi\I$ and~$-1$,
hence the inverse change of variable
$
  q(\ze) \defeq -W_{-1}(- \ee^{\ze})
$
extends analytically to a
neighbourhood of~$\ga_u$, and~\eqref{eqintermed} yields
\beglab{eqsaysnatmaj}
  \frac{z^{z-1}}{\Ga(z)}
  = \frac{1}{2\pi\I} \int_{\ga_u} \ee^{- z \ze}\, q'(\ze)\,\dd\ze
  = \frac{z}{2\pi\I} \int_{\ga_u} \ee^{- z \ze}\, q(\ze)\,\dd\ze
  = \frac{z\,\ee^z}{2\pi\I} \int_{\HKL_0} \ee^{- z \xi}\, q(-1+\xi)\,\dd\xi.
\edla
In the last step, we have used the change of variable $\xi=1+\ze$ with
the convention $\arg(1+\ze(0))=-\pi$, as well as
the fact that the contour $1+\ga_u$ is homotopic to the Hankel
contour~$\HKL_0$ in $\C-2\pi\I\Z$ relatively to $\RE\xi\to+\infty$.

Viewed as a natural-major, $\frac{1}{\I\sqrt{2\pi}}q(-1+\xi)$, which
is none other that~$\chb\chi$,
is thus mapped by~$\cL^0$ to
$\frac{\sqrt{2\pi} z^{z-2}\ee^{-z}}{\Ga(z)}$,
which is $\frac{1}{\la_{-3/2}}$.







Finally, since
$
  \arg\xi = -\pi \iimp
    \chb\chi(\xi) = \frac{\I}{\sqrt{2\pi}} W_{-1}(-\ee^{-1+\xi})
    = -\frac{\I}{\sqrt{2\pi}} q_+(\eipi\xi)
$,
we can use the Puiseux expansion~\eqref{eqPuiseuxWm},
$ 
    q_+(\eipi\xi) = 1 + \sum_{k\ge1} \I^k a_k (2\xi)^{k/2}
    \in\C\{\xi^{1/2}\}$,
    %
  %
  which shows that $\trb\chi \defeq
  \singz{\chb\chi(\xi)}$ is an integrable singularity
  %
  with a minor 
  \[
    \htb\chi(\xi) = -\frac{\I}{\sqrt{2\pi}} \big( q_+(\eipi\xi) -
    q_+(\emipi\xi)\big)
= \frac{\I}{\sqrt{2\pi}} \big( W_{-1}(-\ee^{-1+\xi}) -
    W_0(-\ee^{-1+\xi}) \big) 
    %
    %
    %
\]
whose Puiseux expansion is~(\ref{eqexpchi}\,b).
The proof of Theorem~\ref{thmMinNat} is complete.

\bg

\noindent \textbf{Alien calculus and Stokes phenomenon}   \label{secDirectAlien}

\bg



\noindent A point on the line-segment $(0,\om)$ close
to $\om\in\ti\C$ can be written $\om+\xi$ with the convention
$\arg\xi=\arg\om-\pi$. 
%
%
Given $\trb\ph$, there are finitely many points
$\om_1<\om_2<\cdots<\om_{r-1}$ on $(0,\om)$ such that the analytic continuation of the
minor~$\htb\ph$ exists along every path~$\ga_\eps$,
$\eps\in\{+,-\}^{r-1}$, obtained by following $(0,\om)$ and
circumventing each $\om_j$ to the right (\resp left) if $\eps_j=+$
(\resp $-$) without going backwards. Then \'Ecalle defines two
operators:
%
\begla
\De^+_\om \trb\ph \defeq \singz{\cont_{\ga_{+\cdots+}} \htb\ph(\om+\xi)},
\quad
\De_\om \trb\ph \defeq \sum_\eps \tfrac{p(\eps)!q(\eps)!}{r!}
\singz{\cont_{\ga_\eps} \htb\ph(\om+\xi)}
%
\edla
with $p(\eps)\defeq$ number of~`$+$' in~$\eps$ and $q(\eps)\defeq r-1-p(\eps)$.
The operators~$\De_\om$ are derivations for~$\Dtr *$
and generate an \emph{$\infty$-dimensional free Lie algebra of derivations acting on the algebra of
singularities}, giving rise to interesting algebraic combinatorics
questions~\cite{Eca85}.
The operators~$\De^+_\om$ (which can be expressed in terms of~$\De_{\om'}$, $\om'\in(0,\om)$) satisfy a modified Leibniz
rule and have a natural interpretation in terms of Stokes phenomenon
\cite{Eca81,Eca85,Eca93,MS16},
%
%
which we now illustrate on~$\trn\la_c$. The analytic continuation
of~$W_0$ and~$W_{-1}$ is known well enough so
that~\eqref{eqfirstdefhtnlatdem} allows one to compute directly, with
$\om^\pm_m \defeq 2\pi m\,\ee^{\pm\Itpd}$ for $m\ge1$,
%
\beglab{eqDeptilac}
\De_{\om^\pm_m}\trn\la_c = \pm\tfrac{1}{m}\trn\la_c, \qquad
\De^+_{\om^+_m}\trn\la_c = \trn\la_c, \qquad
\De^+_{\om^-_m}\trn\la_c = \bigg\{
\begin{aligned}
&- \trn\la_c & &\text{for $m = 1$} \\[-.5ex]
&\ens 0  & &\text{for $m \ge2$}
\end{aligned}
%
\edla
with $c=3/2$ and thus also\footnote{A proof can
also be derived from the relation $\ti\la = \exp{\ti\mu}$ as in
\cite{MS16}---beware of the typo in~\cite[(6.99)]{MS16}.} 
with arbitrary~$c$.
It follows that
%
%
\begin{align*}
-\pi<\arg z<0 &\imp
\cL^{J_1} \trn\la_c(z) = \cL^{J_2}\big(
\trn\la_c + \sum_{m\ge1} \ee^{-\om^+_m z} \De^+_{\om^+_m} \trn\la_c \big)
= \tfrac{1}{1-\ee^{-2\pi\I z}}\, \cL^{J_2} \trn\la_c(z) \\[.5ex]
0<\arg z<\pi &\imp
\cL^{J_4} \trn\la_c(z) = \cL^{J_1}\big(
\trn\la_c + \sum_{m\ge1} \ee^{-\om^-_m z} \De^+_{\om^-_m} \trn\la_c \big)
               = (1-\ee^{2\pi\I z}) \cL^{J_1} \trn\la_c(z)
               %
\end{align*}
(which can be viewed as the solution to a Riemann-Hilbert
problem---the one induced by the wall-crossing formula in
Donaldson-Thomas ``Doubled~$A_1$'' theory if $c=0$).
Recall that $\cL^{J_1} \trn\la_{3/2}=\la_{3/2}$.
Since $J_1=-\pi+J_2=\pi+J_4$,
\eqref{eqsymmchilatdem} implies
$\cL^{J_2} \trn\la_{3/2}(z) = -\I\gL^{J_1}\htb\chi(\ee^{\I\pi}z) = -{\I}/{\la_{-3/2}(\ee^{\I\pi}z)}$
and
$\cL^{J_4} \trn\la_{3/2}(z) = \I\gL^{J_1}\htb\chi(\ee^{-\I\pi}z) = {\I}/{\la_{-3/2}(\ee^{-\I\pi}z)}$,
hence both Stokes phenomena amount to the reflection formula
$\Ga(z)\Ga(1-z) = \frac{\pi}{\sin(\pi z)}$.


\section{Proof of Theorem~\ref{thmRealMaj}}   \label{secrealmaj}
\label{secpfthmdeux}

%
%
Consider~$\la(z)$ for $z\in\Rp$.
We know that $\la(z)$ is bounded on $[u,+\infty)$ for any $u>0$ hence,
for any $c\in\C$,
the real-major for $\la_c(z)= z^{-c} \la(z)$ that
is defined by~\eqref{eqrhcbphu} is holomorphic for $\RE\xi>0$ and we
get its analytic continuation to a sector of opening~$3\pi$ in the
Riemann surface of the logarithm~$\ti\C$ by varying the integration half-line.

Let us do it under the assumption $\RE c<\dem$.
Since
$\la_c(z) \underset{z\to0^+}{\sim} \frac{1}{\sqrt{2\pi}} z^{-c-\demi}$
is then integrable at $z=0$,
$u=0$ yields a real-major $\rchn\la_c(\xi) = \gL^0\la_c(\xi)$.
%
%
Moreover,
\beglab{unifestimlac}
\la_c(z) \underset{\abs z\to\infty}{\sim} z^{-c}
\quad\text{and}\quad
\la_c(z) \underset{\abs z\to 0}{\sim} \tfrac{1}{\sqrt{2\pi}} z^{-c-\demi}
\quad
\text{uniformly for} \ens
\abs{\arg z}<\pi-\eps
\edla
for any $0<\eps<\pi$,
hence we can vary the integration ray in the above Laplace transform and
obtain the analytic continuation of~$\rchn\la_c(\xi)$ to a sector of
opening~$3\pi$ in~$\ti\C$ in the form
\begla
\rchn\la_c(\xi) =
\gL^{J_1}\la_c(\xi)
\quad\text{for}\;
-\ttrpd < \arg\xi < \ttrpd,
\quad\text{with} \; J_1 \defeq \tintI.
\edla
Moreover, in view of~\eqref{unifestimlac}, this representation of $\rchn\la_c(\xi)$ shows that it is
bounded in the domain $\{ \abs\xi >\tau \;\text{and}\; \abs{\arg\xi} <
\trpd-\eps \}$ for every $\tau,\eps>0$.
Therefore Fourier-Laplace inversion yields
%
%
$\la_c(z) = \rchnLth\rchn\la_c(z)$
%
  for all $z\in\C-\R_{\le0}$ 
  and $\tht\in J_1$ 
  such that 
  $\tht+\arg z \in J_1$. 
%
%


%
We now prove Formula~\eqref{eqrchnlac}.
%
%
For each $z>0$,
\[
  \Ga(z) = \int_0^{+\infty} \ee^{-x} x^z \,\tfrac{\dd
    x}{x} = z^z \int_\R \ee^{-z\ee^Q} \ee^{zQ}\,\dd Q
\]
(change of variable $x=z\,\ee^Q$),
hence 
$\dst 
\la_c(z) = \tfrac{1}{\sqrt{2\pi}} z^{-c+\demi}\int_\R
\ee^{-z(\ee^Q-Q-1)} \,\dd Q
%
%
$ and
\[
\rchn\la_c(\xi) = \tfrac{1}{\sqrt{2\pi}} 
\int_0^{+\infty} \!\dd z \int_\R \!\dd Q\, 
z^{-c+\demi} \, \ee^{-z(\xi+\ee^Q-Q-1)}.
\]
Now $\RE(\xi+\ee^Q-Q-1)>0$ by~\eqref{ineqexp} and $\RE(\tdemi-c)>0$, thus 
\[ 
\int_0^{+\infty} \! z^{-c+\demi} \, \ee^{-z(\xi+\ee^Q-Q-1)}\dd z
= \Ga(\tdem-c)(\xi+\ee^Q-Q-1)^{c-\tdemi}
\]
by~\eqref{eqLaplMon}, and the conclusion follows.
%


%
To conclude the proof of Theorem~\ref{thmRealMaj}, we observe that,
to follow the analytic continuation of~$\rchn\la_c(\xi)$ along a
path~$\Ga$ of~$\ti\C$ that starts on $\{\arg\xi=0\}$, it is
sufficient to deform continuously the integration path
in~\eqref{eqrchnlac} so that the integration variable~$Q$ avoids the
zeroes of the expression $\xi+\ee^Q - Q- 1$.
Such a continuous deformation exists as long as~$0$ is not a critical
value of $Q \mapsto \xi+\ee^Q - Q- 1$,
\ie as long as the path~$\Ga$ avoids any critical value of $Q \mapsto
-\ee^Q + Q + 1$.
These critical values are precisely the points of $2\pi\I\Z$.

\bg

\bg

\subsubsection*{Acknowledgements}

The author thanks Capital Normal University for their hospitality
during the period September 2019--February 2020, where this research
was started.
%
This paper is partly a result of the ERC-SyG project, Recursive and
Exact New Quantum Theory (ReNewQuantum) which received funding from
the European Research Council (ERC) under the European Union's Horizon
2020 research and innovation programme under grant agreement No
810573.


\newpage

















\printbibliography


\end{document}